\documentclass[12pt]{amsart}
\usepackage[margin=1in]{geometry}

\usepackage{mathtools, thmtools}
\usepackage{amssymb, mleftright, bm}
\usepackage[backend=biber, style=numeric, sorting=nyt, giveninits=true]{biblatex}
\usepackage[stretch=10, shrink=10]{microtype}
\usepackage{graphicx, subcaption}
\usepackage[bb=dsserif, bbscaled=1.1]{mathalpha}
\usepackage[T1]{fontenc}
\usepackage{lmodern}
\usepackage{hyperref}

\DeclareMathOperator{\pv}{p.v.}

\NewDocumentCommand{\abs}{m}{\mleft\lvert #1 \mright\rvert}
\NewDocumentCommand{\df}{}{\mathop{}\!\mathrm{d}}

\RenewDocumentCommand{\Re}{}{\operatorname{Re}}
\RenewDocumentCommand{\Im}{}{\operatorname{Im}}

\numberwithin{equation}{section}
\declaretheorem[numberwithin=section]{theorem}
\declaretheorem[numberlike=theorem]{lemma, corollary, proposition}

\graphicspath{ {./image} }

\begin{document}

\title{Discrete Fourier Transform and~\( L \)-functions}
\author{Di Liu}
\date{}

\begin{abstract}
  We give the converse to Dirichlet's theorem on primes
  in arithmetic progressions by generalizing an old result
  of~\citeauthor{guinandSummationFormulaTheory1948}.
\end{abstract}

\maketitle

\section{Introduction}

In~\citeyear{riemannUeber1859}, \citeauthor{riemannUeber1859}~\cite{riemannUeber1859}
revealed the Fourier duality between prime numbers
and nontrivial zeros of~\( \zeta\mleft( s \mright) \), denoted by~\( \rho = \beta + i \gamma \).
The \emph{Riemann~Hypothesis}~(RH) asserts~\( \beta = 1 / 2 \).
More explicitly they are related by the following pair of formulae
\begin{align}
  N\mleft( T \mright)
  & \coloneqq \sum_{0 < \Im \rho \le T} 1
  = \frac{1}{\pi} \arg \Gamma\mleft( \frac{1 / 2 + i T}{2} \mright)
  - \frac{\log \pi}{2 \pi} T + 1
  + \arg \zeta\mleft( 1 / 2 + i T \mright),
  \label{nt} \\
  \psi\mleft( x \mright)
  & \coloneqq \sum_{n \le x} \Lambda\mleft( n \mright)
  = x - \frac{\zeta'}{\zeta}\mleft( 0 \mright)
  - \log\mleft( 1 - x^{- 2} \mright)
  - \sum_{\rho} \frac{x^{\rho}}{\rho},
  \label{psi}
\end{align}
where~\( \Lambda\mleft( n \mright) \) is the von~Mangoldt function,
and halving convention is adopted whenever the sum cutoff coincides with a term.
In~\eqref{nt} the primes are hidden in the last argument.

Next in~\citeyear{landauUeberNullstellenZetafunktion1912},
\citeauthor{landauUeberNullstellenZetafunktion1912}~\cite{landauUeberNullstellenZetafunktion1912}
showed unconditionally for any fixed~\( x > 1 \) as~\( T \to + \infty \),
\begin{equation}
  \label{ld}
  \sum_{0 < \rho \le T} x^{\rho} = - \frac{T}{2 \pi} \Lambda\mleft( x \mright)
  + O\mleft( \log T \mright).
\end{equation}
This says zeros are able to locate prime powers in addition to
the counting in~\eqref{psi}.
Analytically it approximates~\( \psi'\mleft( x \mright) \),
as spikes of size~\( T \) behave like Dirac deltas.

Later in~\citeyear{inghamDistributionPrimeNumbers2000},
\citeauthor{inghamDistributionPrimeNumbers2000}~\cite[81]{inghamDistributionPrimeNumbers2000}
suggested comparing~\eqref{psi} to the counting function
\begin{equation}
  \label{st}
  \sum_{n \le x} 1 = x - 1 / 2
  + \frac{1}{\pi} \sum_{n = 1}^{\infty} \frac{\sin\mleft( 2 \pi n x \mright)}{n},
\end{equation}
a modified version of the sawtooth wave.
In~\citeyear{guinandPoissonsSummationFormula1941},
\citeauthor{guinandPoissonsSummationFormula1941}~\cite{guinandPoissonsSummationFormula1941}
proved a general version of the Poisson summation formula,
equivalent to the derivative of~\eqref{st}.

Then in~\citeyear{guinandSummationFormulaTheory1948},
under RH
\citeauthor{guinandSummationFormulaTheory1948}~\cite[Thm.~4]{guinandSummationFormulaTheory1948}
showed for any fixed~\( y \in \mathbb{R} \)
as~\( x \to + \infty \),
\begin{multline}
  \label{gnd}
  \sum_{0 < n \log p < x}
  \frac{\log p}{p^{n / 2}} \cos\mleft( y n \log p \mright)
  - \frac{e^{x / 2}}{1 / 4 + y^{2}}
  \mleft( \tfrac{1}{2} \cos\mleft( y x \mright) + y \sin \mleft( y x \mright) \mright) \\
  - \biggl( \, \sum_{0 < n \log p < x}
  \frac{\log p}{p^{n / 2}} - 2 e^{x / 2} \biggr)
  \cos\mleft( y x \mright)
  =
  \begin{cases}
    - x + O\mleft( 1 \mright) & \text{if}~y = \gamma, \\
    O\mleft( 1 \mright) & \text{otherwise}.
  \end{cases}
\end{multline}
This remarkable theorem is essentially
the zero indicator function without any smoothing,
and the dual to Landau's result~\eqref{ld}
as observed by~\citeauthor{guinandSummationFormulaTheory1948}.
Analytically it approximates~\( N'\mleft( T \mright) \) hence relates to
the derivative of~\( S\mleft( T \mright) \coloneqq
\arg \zeta\mleft( 1 / 2 + i T \mright) \) via~\eqref{nt}.
It also shows that the distributional Fourier transform
\begin{equation}
  \label{dstFt}
  \sum_{\gamma} \delta\mleft( t - \gamma \mright)
  - \Re \frac{\Gamma'}{\Gamma}\mleft( \frac{1 / 2 + i t}{2} \mright)
  + \log \pi
  \leftrightarrow
  - \sum_{p} \sum_{n \in \mathbb{Z} \setminus \mleft\{ 0 \mright\}}
  \frac{\log p}{p^{n / 2}} \cdot \delta\mleft( x - n \log p \mright)
  + 2 \cosh\mleft( \frac{x}{2} \mright)
\end{equation}
implicit in the explicit formula in the same paper and
\citeauthor{weilFormulesExplicitesTheorie1952}~\cite{weilFormulesExplicitesTheorie1952}
almost holds true as is.
Such transforms are called \emph{Concordance}
by~\citeauthor{guinandConcordanceHarmonicAnalysis1959}~\cite{guinandConcordanceHarmonicAnalysis1959}
and more recently, \emph{Crystalline~Measure}
by~\citeauthor{meyerMeasuresLocallyFinite2016}~\cite{meyerMeasuresLocallyFinite2016}.
In fact, the functions in~\eqref{dstFt} are not quite crystalline
as they are not pure Dirac deltas. However, the classical example
~\( \sum_{n \in \mathbb{Z}} \delta\mleft( x - n \mright) \) being self dual
is equivalent to~\( \sum_{n \in \mathbb{Z} \setminus \mleft\{ 0 \mright\}}
\delta\mleft( x - n \mright) - 1 \) being so as well.
Therefore the~\( 0 \)-component, also known as
the \emph{Direct~Current}~(DC) component,
could be regarded as background average.

A drawback of~\citeauthor{guinandConcordanceHarmonicAnalysis1959}'s
result~\eqref{gnd} is the oscillating sum of size~\( O\mleft( x^{2} \mright) \)
on the second line by~\cite[(2.14)]{guinandSummationFormulaTheory1948},
larger than the spikes of size~\( x \) at zeros.
In this paper alternative versions are proved
for the Riemann zeta function and Dirichlet~\( L \)-functions.
Nontrivial zeros of an~\( L\mleft( s, \chi \mright) \)
are denoted by~\( \rho_{\chi} = \beta_{\chi} + i \gamma_{\chi} \),
and the \emph{Generalized~Riemann~Hypothesis}~(GRH)
asserts~\( \beta_{\chi} = 1 / 2 \).
Euler's totient function is denoted by~\( \varphi\mleft( n \mright) \).

\begin{theorem}
  \label{tmZt}
  Assume RH, let~\( \delta_{\zeta}\mleft( y \mright) = k_{\gamma} \)
  if~\( y = \gamma \) for a zero~\( 1 / 2 + i \gamma \)
  of order~\( k_{\gamma} \), and~\( 0 \) otherwise.
  Let~\( x > e \) and~\( \varepsilon_{x} = 3 \log x / x \),
  then for~\( y \in \mathbb{R} \) as~\( x \to + \infty \),
  \begin{multline}
    \label{indZt}
    \Re \mleft( \sum_{0 < n \log p \le x}
      \frac{\log p}{p^{n \mleft( 1 / 2 + \varepsilon_{x} + i y \mright)}}
      - \frac{e^{x \mleft( 1 / 2 - \varepsilon_{x} - i y \mright)}}{1 / 2 - \varepsilon_{x} - i y} \mright)
    - \frac{1}{2} \Re \frac{\Gamma'}{\Gamma}\mleft( \frac{1 / 2 - i y}{2} \mright)
    + \frac{1}{2} \log \pi \\
    = - \frac{x}{3 \log x} \cdot \delta_{\zeta}\mleft( y \mright)
    + O\mleft( \frac{\log x}{x} \cdot
    \frac{\log \mleft( 2 + \abs{ y } \mright)}
    {\min_{\gamma \ne y} \abs{ y - \gamma }^{2}}\mright).
  \end{multline}
\end{theorem}

\begin{theorem}
  \label{tmDr}
  Let~\( \chi \bmod q \) be primitive,
  \( \mathfrak{a}_{\chi} \) its parity according
  to~\( \chi\mleft( -1 \mright) = \mleft( - 1 \mright)^{\mathfrak{a}_{\chi}} \).
  Assume~GRH for~\( L\mleft( s, \chi \mright) \),
  let~\( \delta_{\chi}\mleft( y \mright) = k_{\gamma_{\chi}} \)
  if~\( y = \gamma_{\chi} \) for a zero~\( 1 / 2 + i \gamma_{\chi} \)
  of order~\( k_{\gamma_{\chi}} \), and~\( 0 \) otherwise.
  Let~\( x > e \) and~\( \varepsilon_{x} = 3 \log x / x \),
  then for~\( y \in \mathbb{R} \) as~\( x \to + \infty \),
  \begin{multline}
    \label{indDr}
    \Re \sum_{0 < n \log p \le x}
    \frac{\chi\mleft( p^{n} \mright) \log p}{p^{n \mleft( 1 / 2 + \varepsilon_{x} + i y \mright)}}
    - \frac{1}{2} \Re \frac{\Gamma'}{\Gamma}\mleft( \frac{1 / 2 + \mathfrak{a}_{\chi} - i y}{2} \mright)
    + \frac{1}{2} \log \frac{\pi}{q} \\
    = - \frac{x}{3 \log x} \cdot \delta_{\chi}\mleft( y \mright)
    + O_{q}\mleft( \frac{\log x}{x} \cdot
      \frac{\log \mleft( 2 + \abs{ y } \mright)}
      {\min_{\gamma_{\chi} \ne y} \abs{ y - \gamma_{\chi} }^{2}} \mright).
  \end{multline}
\end{theorem}

The downside here is the less than ideal spike size.
On the other hand, the dependence on~\( y \) is made explicit
and cancels with the smooth terms of~\( N'\mleft( T \mright) \)
in~\eqref{nt}, giving a cleaner approximation to~\( S'\mleft( T \mright) \).
Moreover, by moving off the critical line,
the origin is treated uniformly instead of being
excluded in~\citeauthor{landauUeberNullstellenZetafunktion1912}'s~\eqref{ld},
or subtracted in~\citeauthor{guinandSummationFormulaTheory1948}'s~\eqref{gnd}.
At~\( y = 0 \) the trignometric series~\eqref{indZt}
essentially reduces to the sum~\eqref{psi} with a different weight,
which further demonstrates the intertwining nature between
those counting formulae and their indicator versions.
Usually indicators are regarded as derivatives of counting formulae,
yet counting are indicators of their dual evaluated at the origin.

A corollary immediately follows by some linear algebra.

\begin{corollary}
  \label{crInd}
  Let~\( q \ge 3 \) and let~\( \mleft\{ \chi \mright\} \) run through the~\( \varphi\mleft( q \mright) \)
  Dirichlet characters modulo~\( q \),
  \( \mleft\{ \mathfrak{a}_{\chi} \mright\} \) their parity according
  to~\( \chi\mleft( -1 \mright) = \mleft( -1 \mright)^{\mathfrak{a}_{\chi}} \) respectively.
  Assume~GRH for those~\( L\mleft( s, \chi \mright) \)'s,
  let~\( \delta_{q}\mleft( y \mright) = \sum_{\chi} k_{\gamma_{\chi}} \)
  if~\( y = \gamma_{\chi} \) for a zero~\( 1 / 2 + i \gamma_{\chi} \)
  of order~\( k_{\gamma_{\chi}} \) of~\( L\mleft( s, \chi \mright) \),
  and~\( 0 \) otherwise.
  Let \( x > e \)~and~\( \varepsilon_{x} = 3 \log x / x \),
  then for~\( y \in \mathbb{R} \) as~\( x \to + \infty \),
  \begin{multline}
    \label{indDec}
    \Re \Biggl( \: \sum_{\substack{ 0 < n \log p \le x \\ p^{n} \equiv 1 \bmod q }}
      \frac{\log p}{p^{n \mleft( 1 / 2 + \varepsilon_{x} + i y \mright)}}
      - \frac{1}{\varphi\mleft( q \mright)} \cdot
      \frac{e^{x \mleft( 1 / 2 - \varepsilon_{x} - i y \mright)}}
    {1 / 2 - \varepsilon_{x} - i y} \Biggr) \\
    \shoveright{ - \frac{1}{4} \Re \frac{\Gamma'}{\Gamma}\mleft( \frac{1 / 2 - i y}{2} \mright)
      - \frac{1}{4} \Re \frac{\Gamma'}{\Gamma}\mleft( \frac{3 / 2 - i y}{2} \mright)
    + \frac{1}{2} \log \frac{\pi}{q} } \\
    =
    - \frac{1}{\varphi\mleft( q \mright)} \cdot \frac{x}{3 \log x}
    \cdot \delta_{q}\mleft( y \mright) + O_{q}\biggl( 1 + \frac{\log x}{x} \cdot
      \frac{\log \mleft( 2 + \abs{ y } \mright)}
      {\min_{\substack{ \gamma_{\chi} \ne y \\ \chi \bmod q }}
    \abs{ y - \gamma_{\chi} }^{2}} \biggr).
  \end{multline}
  In particular, zeros of~\( L\mleft( s, \chi \mright) \)
  for~\( \chi \) the principal character coincide with
  those of the Riemann zeta function.
\end{corollary}

This is the converse to Dirichlet's theorem
on primes in arithmetic progressions,
or the zero indicator for the Dedekind~zeta~function~\( \zeta_{K}\mleft( s \mright) \) of a cyclotomic field.
It says that counterintuitively the \emph{fewer} primes used,
the \emph{more} zeros detected.
As will be shown, this is not unique
to the residue class~\( 1 \bmod q \) but all of those coprime to~\( q \).

The proof uses the \emph{character table} of the
group~\( \mleft( \mathbb{Z} / q \mathbb{Z} \mright)^{\times} \),
while its additive version is more widely known as
the \emph{Discrete~Fourier~Transform}~(DFT)~matrix,
a fundamental object in both engineering and mathematics
as discussed in~\citeauthor{auslanderComputingFiniteFourier1979}~\cite{auslanderComputingFiniteFourier1979}.
While individual~\( L \)-functions are transforms of prime powers
via~\eqref{dstFt} and its variants, the DFT of them in residue classes by \emph{twisting}
are functions determined by zeros of~\( L \)-functions in some family,
with the Riemann zeta function always being the 0-component.

From a signal processing point of view,
\eqref{indDec} describes the DC output of
a \( \varphi\mleft( q \mright) \)-point DFT system.
The basic law governing such systems is the sampling theorem
due to \citeauthor{whittakerXVIIIonFunctionsWhich1915}~\cite{whittakerXVIIIonFunctionsWhich1915}
and \citeauthor{shannonCommunicationPresenceNoise1949}~\cite{shannonCommunicationPresenceNoise1949},
saying that for a bandlimited signal, once the frequency~\( \varphi\mleft( q \mright) \)
exceeds its Nyquist rate~\cite{nyquistCertainTopicsTelegraph1928}
the output should stablize, which is not the case here.
Consequently this implies our target has infinite bandwidth, always \emph{undersampled},
therefore the Riemann zeta function as well as all Dirichlet~\( L \)-functions
arise as its \emph{aliasing}. This agrees with the uniform treatment
in~\citeauthor{bohr1913}~\cite{bohr1913} analytically,
as well as in~\citeauthor{tateFourierAnalysisNumber1950}'s
thesis~\cite{tateFourierAnalysisNumber1950} and Iwasawa theory algebraically.
It also explains how non-primitive~Dirichlet \( L \)-functions
are induced by primitive ones in a similar fashion,
and the extra error~\( O_{q}\mleft( 1 \mright) \) is due to their difference.

\section{Proof}

The first lemma needed is Corollary~(b)
in~\citeauthor{davenportMultiplicativeNumberTheory1980}~\cite[99]{davenportMultiplicativeNumberTheory1980}
where we use the bound~\( \log \mleft( 2 + \abs{ y } \mright) \) instead of~\( \log \abs{ y } \),
as digamma is analytic on~\( \Re \mleft( s \mright) = 2 \).
\begin{lemma}
  \label{lmZ}
  If~\( \rho = \beta + i \gamma \) runs through the nontrivial zeros
  of~\( \zeta\mleft( s \mright) \), then for~\( y \in \mathbb{R} \)
  \begin{equation*}
    \sum_{\abs{ \Im \rho - y } \ge 1} \frac{1}{\mleft( y - \gamma \mright)^{2}}
    = O\mleft( \log \mleft( 2 + \abs{ y } \mright) \mright).
  \end{equation*}
\end{lemma}
The second is the version for Dirichlet~\( L \)-functions,
following from the Lemma
in~\cite[102]{davenportMultiplicativeNumberTheory1980}
by the same argument.
\begin{lemma}
  \label{lmC}
  If~\( \rho_{\chi} = \beta_{\chi} + i \gamma_{\chi} \) runs through the nontrivial zeros
  of~\( L\mleft( s, \chi \mright) \), where~\( \chi \bmod q \) is primitive,
  then for~\( y \in \mathbb{R} \)
  \begin{equation*}
    \sum_{\abs{ \Im \rho_{\chi} - y } \ge 1} \frac{1}{\mleft( y - \gamma_{\chi} \mright)^{2}}
    = O_{q}\mleft( \log \mleft( 2 + \abs{ y } \mright) \mright).
  \end{equation*}
\end{lemma}
The third is the Lemma in~\cite[105]{davenportMultiplicativeNumberTheory1980}.
\begin{lemma}
  \label{lmT}
  Denote the Heaviside step function by~\( H\mleft( u \mright) \), and let
  \begin{equation*}
    I\mleft( u, T \mright) = \frac{1}{2 \pi i}
    \int_{c - i T}^{c + i T} \frac{u^{s}}{s} \df s.
  \end{equation*}
  Then for~\( u > 0, c > 0, T > 0 \),
  \begin{equation*}
    \abs{ I\mleft( u, T \mright) - H\mleft( u - 1 \mright) } \le
    \begin{cases}
      u^{c} \min \mleft\{ 1, \mleft( T \abs{ \log u } \mright)^{- 1} \mright\}
      & \text{if}~u \ne 1, \\
      c / T & \text{otherwise}.
    \end{cases}
  \end{equation*}
\end{lemma}
The last one is~\citeauthor{inghamDistributionPrimeNumbers2000}~\cite[Thm.~26]{inghamDistributionPrimeNumbers2000}, where we treat the positive and negative intervals simultaneously.
In the worst case the number of zeros are doubled in the interval~\( \mleft( n, n + 1 \mright] \),
giving an extra constant~\( 2 \) absorbed in~\( \ll \).
\begin{lemma}
  \label{lmH}
  For~\( - 1 \le \sigma \le 2 \) there exists a sequence of
  heights~\( \mleft\{ T_{t, n} \mright\} \) for~\( n \ge 2 \) such that
  \begin{equation*}
    \abs{ \frac{\zeta'}{\zeta}\mleft( \sigma + i t \pm i T_{t, n} \mright) }
    \ll \log^{2} \mleft( \abs{ t } + T_{t, n} \mright)
    ~\text{and}~
    \abs{ t } + n < T_{t, n} < \abs{ t } + n + 1.
  \end{equation*}
\end{lemma}

\subsection{Proof~of~\autoref{tmZt}}

We start with the characteristic function of
the interval~\( \mleft[ 1, x \mright] \) and its Mellin transform,
\begin{equation}
  \label{mel}
  \chi_{\mleft[ 1, x \mright]}\mleft( u \mright)
  \leftrightarrow
  w_{x}\mleft( s \mright) \coloneqq
  \begin{dcases}
    \log x & \text{if}~s = 0, \\
    \frac{x^{s} - 1}{s} & \text{otherwise}.
  \end{dcases}
\end{equation}
Consider the following Dirichlet polynomial
defined inside the strip~\( \mleft\{ z \in \mathbb{C} \colon 0 \le \Im z < 1 / 2 \mright\} \),
\begin{equation}
  \label{dpDef}
  D_{x}\mleft( z \mright)
  \coloneqq \sum_{n \le x}
  \frac{\Lambda\mleft( n \mright)}{n^{1 / 2 - i z}}
  = \frac{1}{2 \pi i} \int_{\mleft( c \mright)}
  - \frac{\zeta'}{\zeta}\mleft( 1 / 2 - i z + s \mright)
  \cdot w_{x}\mleft( s \mright) \df s,
\end{equation}
where \( c = 1 + 1 / \log x \) and we write~\( z = - y + i \varepsilon \).
The integral formula follows by Mellin inversion applied to the pair~\eqref{mel},
and we are interested in its real part.
We choose a~\( T \) such that~\( T / \log^{2} T > \abs{ y } \)~and \( y \pm T \) not a zero,
then residue theorem says
\begin{multline}
  \label{rsd}
  \frac{1}{2 \pi i}
  \mleft( \int_{c - i T}^{c + i T} + \int_{c + i T}^{- c + i T}
  + \int_{- c + i T}^{- c - i T} + \int_{- c - i T}^{c - i T} \mright)
  - \frac{\zeta'}{\zeta}\mleft( 1 / 2 - i z + s \mright)
  \cdot w_{x}\mleft( s \mright) \df s \\
  = w_{x}\mleft( 1 / 2 + i z \mright)
  - \sum_{\abs{ \Im \rho - y } < T} w_{x}\mleft( \rho - 1 / 2 + iz \mright),
\end{multline}
and the integrals are denoted~\( I_{1}, \ldots, I_{4} \) in the order of appearance.
For the two horizontal parts \( I_{2} \)~and~\( I_{4} \),
as~\( T \) runs through the sequence in~\autoref{lmH},
\begin{equation*}
  I_{2} + I_{4}
  \ll \frac{\log^{2}\mleft( \abs{ y } + T \mright)}{T}
  \int_{- c}^{c} x^{\sigma} + 1 \df \sigma
  \ll \frac{x \log^{2} T}{T \log x}.
\end{equation*}
For~\( I_{3} \), the functional equation
\begin{equation}
  \label{fneq}
  - \frac{\zeta'}{\zeta}\mleft( 1 / 2 + s \mright)
  - \frac{\zeta'}{\zeta}\mleft( 1 / 2 - s \mright)
  = \frac{1}{2} \frac{\Gamma'}{\Gamma}\mleft( \frac{1 / 2 + s}{2} \mright)
  + \frac{1}{2} \frac{\Gamma'}{\Gamma}\mleft( \frac{1 / 2 - s}{2} \mright)
  - \log \pi,
\end{equation}
and an \( s \mapsto - s \)~change of variable gives
\begin{multline*}
  I_{3} = \frac{1}{2 \pi i} \int_{c - i T}^{c + i T}
  - \frac{\zeta'}{\zeta}\mleft( 1 / 2 + i z + s \mright)
  \cdot w_{x}\mleft( - s \mright) \df s \\
  - \frac{1}{2 \pi i} \int_{- c - i T}^{- c + i T}
  \mleft( \frac{1}{2} \frac{\Gamma'}{\Gamma}\mleft( \frac{1 / 2 - i z + s}{2} \mright)
  + \frac{1}{2} \frac{\Gamma'}{\Gamma}\mleft( \frac{1 / 2 + i z - s}{2} \mright)
  - \log \pi \mright) \cdot w_{x}\mleft( s \mright) \df s.
\end{multline*}
We move the digamma integral to~\( \Re \mleft( s \mright) = 0 \),
picking up the pole at~\( s = - 1 / 2 + i z \),
and note that~\( \abs{ \Gamma'/\Gamma\mleft( \sigma + i t \mright) } \ll \log t \)
in any fixed vertical strip.
With the identity~\( w_{x}\mleft( it \mright) + w_{x}\mleft( - it \mright)
= 2 \sin\mleft( t \log x \mright) / t \), the above reduces to
\begin{multline}
  \label{sinc}
  I_{3} = \frac{1}{2 \pi i} \int_{c - i T}^{c + i T}
  - \frac{\zeta'}{\zeta}\mleft( 1 / 2 + i z + s \mright)
  \cdot w_{x}\mleft( - s \mright) \df s
  - w_{x}\mleft( - 1 / 2 + i z \mright)
  + O\mleft( \frac{\log T}{T} \mright) \\
  - \frac{1}{2 \pi} \int_{- T}^{T}
  \mleft( \frac{1}{2}\frac{\Gamma'}{\Gamma}
    \mleft( \frac{1 / 2 + i z + it}{2} \mright)
  + \frac{1}{2} \frac{\Gamma'}{\Gamma}
  \mleft( \frac{1 / 2 - i z - it}{2} \mright)
  - \log \pi \mright) \cdot \frac{\sin\mleft( t \log x \mright)}{t} \df t.
\end{multline}
Together the residue formula~\eqref{rsd} becomes
\begin{multline}
  \label{rsd13}
  \frac{1}{2 \pi i} \int_{c - i T}^{c + i T}
  - \frac{\zeta'}{\zeta}\mleft( 1/2 - i z + s \mright)
  \cdot w_{x}\mleft( s \mright)
  - \frac{\zeta'}{\zeta}\mleft( 1/2 + i z + s \mright)
  \cdot w_{x}\mleft( -s \mright) \df s \\
  \shoveleft = w_{x}\mleft( 1 / 2 + i z \mright)
  + w_{x}\mleft( - 1 / 2 + i z \mright)
  - \sum_{\abs{ \Im \rho - y } < T} w_{x}\mleft( \rho - 1 / 2 + i z \mright)
  + O\mleft( \frac{x \log^{2} T}{T \log x} \mright) \\
  + \frac{1}{2 \pi} \int_{- T}^{T}
  \mleft( \frac{1}{2} \frac{\Gamma'}{\Gamma}
    \mleft( \frac{1 / 2 + i z + it}{2} \mright)
  + \frac{1}{2} \frac{\Gamma'}{\Gamma}
  \mleft( \frac{1 / 2 - i z - it}{2} \mright)
  - \log \pi \mright) \cdot \frac{\sin\mleft( t \log x \mright)}{t} \df t.
\end{multline}

Now we estimate the difference between the two integrals on the left
and the Dirichlet polynomial~\( D_{x}\mleft( z \mright) \), that is,
the error introduced by truncating the Mellin inversion~\eqref{dpDef} at height~\( T \).
This part follows~\citeauthor{davenportMultiplicativeNumberTheory1980}~\cite[106-107]{davenportMultiplicativeNumberTheory1980}.
In terms of~\( \zeta'/\zeta \), the lines of integraton
are~\( \Re \mleft( 1 / 2 \pm iz + s \mright)
  = 3 / 2 \pm \varepsilon + 1 / \log x > 1 + 1 / \log x \),
therefore by absolute convergence they equal
\begin{equation*}
  \sum_{n \ge 2} \frac{\Lambda\mleft( n \mright)}{n^{1 / 2 - i z}}
  \cdot \frac{1}{2 \pi i} \int_{c - i T}^{c + i T}
  \frac{\mleft( x / n \mright)^{s} - n^{-s}}{s} \df s
  + \sum_{m \ge 2} \frac{\Lambda\mleft( m \mright)}{m^{1 / 2 + i z}}
  \cdot \frac{1}{2 \pi i} \int_{c - i T}^{c + i T}
  \frac{m^{- s} - \mleft( m x \mright)^{- s}}{s} \df s.
\end{equation*}
By~\autoref{lmT}, its difference from summing with
the Heaviside step function is at most
\begin{multline*}
  \sum_{n \ne x} \frac{\Lambda\mleft( n \mright)}{n^{1 / 2 + \varepsilon}}
  \mleft( \mleft( \frac{x}{n} \mright)^{c}
  \min \mleft\{ 1, \frac{1}{T \abs{ \log \mleft( x / n \mright) }} \mright\}
  + \frac{1}{T n^{c} \log n} \mright)
  + \frac{1 + \log x}{T x^{1 / 2 + \varepsilon}}
  + \frac{1}{e T x^{3/2 + \varepsilon}} \\
  + \sum_{m \ge 2} \frac{\Lambda\mleft( m \mright)}{m^{1 / 2 - \varepsilon}}
  \mleft( \frac{1}{T m^{c} \log m}
    + \frac{1}{T \mleft( m x \mright)^{c} \log \mleft( m x \mright)} \mright),
\end{multline*}
which reduces to
\begin{equation*}
  \sum_{n \ne x} \frac{\Lambda\mleft( n \mright)}{n^{1 / 2 + \varepsilon}}
  \mleft( \mleft( \frac{x}{n} \mright)^{c} \min \mleft\{ 1,
  \frac{1}{T \abs{ \log \mleft( x / n \mright) }} \mright\} \mright)
  + O\mleft( \frac{\log x}{T} \mright).
\end{equation*}
The contribution from those~\( n \) away from~\( x \),
say~\( \abs{ n / x - 1 } \ge 1 / 4 \), is at most
\begin{equation*}
  \frac{x^{c}}{T \log\mleft( 5 / 4 \mright)}
  \sum_{\substack{ n \ge 2 \\ \abs{ n / x - 1 } \ge 1/4 }}
  \frac{\Lambda\mleft( n \mright)}{n^{3/2 + \varepsilon}}
  \ll \frac{x}{T}.
\end{equation*}
For those close to and smaller than~\( x \),
say~\( 3 / 4 < n / x < 1 \),
let~\( x_{1} \) be the largest prime power strictly smaller than~\( x \).
By~\( \log \mleft( x / x_{1} \mright) \ge \mleft( x - x_{1} \mright) / x \),
this single term gives at most
\begin{equation*}
  \frac{\Lambda\mleft( x_{1} \mright)}{x_{1}^{1 / 2 + \varepsilon}}
  \mleft( \frac{4}{3} \mright)^{c}
  \min \mleft\{ 1, \frac{x}{T\mleft( x - x_{1} \mright)} \mright\}
  \ll \frac{\log x}{x^{1 / 2 + \varepsilon}}.
\end{equation*}
For~\( \tfrac{3}{4} x < n < x_{1} \),
using~\( \log \mleft( x / n \mright) \ge \mleft( x_{1} - n \mright) / x_{1} \),
they are bounded by
\begin{equation*}
  \mleft( \frac{4}{3} \mright)^{c}
  \sum_{\tfrac{3}{4} x < n < x_{1}}
  \frac{\Lambda\mleft( n \mright)}{n^{1 / 2 + \varepsilon}}
  \cdot \frac{x_{1}}{T \mleft( x_{1} - n \mright)}
  \ll \frac{x \log x}{T}
  \sum_{1 < m < x / 4} \frac{1}{m \mleft( \frac{3}{4} x + m \mright)^{1/2 + \varepsilon}}
  \ll \frac{x^{1 / 2 - \varepsilon} \log^{2} x}{T}.
\end{equation*}
The symmetric part~\( 1 < n / x < 5/4 \) is similarly treated.
Then provided~\( T \ll x^{3 / 2} \), \eqref{rsd13} reduces to
\begin{multline}
  \label{ztSum}
  \sum_{n \le x} \frac{\Lambda\mleft( n \mright)}{n^{1 / 2 - i z}}
  = w_{x}\mleft( 1 / 2 + i z \mright)
  + w_{x}\mleft( - 1 / 2 + i z \mright)
  - \sum_{\abs{ \Im \rho - y } < T} w_{x}\mleft( \rho - 1 / 2 + i z \mright)
  + O\mleft( \frac{x \log^{2} \mleft( x T \mright)}{T \log x} \mright) \\
  + \frac{1}{2 \pi} \int_{- T}^{T}
  \mleft( \frac{1}{2} \frac{\Gamma'}{\Gamma}
    \mleft( \frac{1 / 2 + i z + it}{2} \mright)
  + \frac{1}{2} \frac{\Gamma'}{\Gamma}
  \mleft( \frac{1 / 2 - i z - it}{2} \mright)
  - \log \pi \mright) \cdot \frac{\sin\mleft( t \log x \mright)}{t} \df t.
\end{multline}
It follows from the step condition in~\autoref{lmH} that
around any large height~\( T \) for all~\( z \)
satisfying~\( T / \log^{2} T > \abs{ y } \),
one can find an index~\( n_{y} \)
such that~\( \abs{ T - T_{y, n_{y}} } \le 1 \).
Therefore changing the range~\( \abs{ \Im \rho - y } < T_{y, n_{y}} \)
to~\( \abs{ \Im \rho } < T \) amounts to adding or subtracting
zeros in the range~\( \abs{ \Im \rho \pm T } \le \abs{ y } + 1  \),
which are far away from~\( z \).
The number of such zeros is~\( \ll \mleft( \abs{ y } + 1 \mright)
\cdot \log \mleft( T + \abs{ y } \mright) \),
and their contribution is at most
\begin{equation*}
  \sum_{\abs{ \Im \rho \pm T } \le \abs{ y } + 1}
  \abs{ w_{x}\mleft( \rho - 1 / 2 + i z \mright) }
  \ll \frac{\mleft( \abs{ y } + 1 \mright) \log \mleft( T + \abs{ y } \mright)}
  {\abs{ T \pm 2 \abs{ y } \pm 1 }}
  \ll \frac{1}{\log T},
\end{equation*}
whence in~\eqref{ztSum} the zero sum range
could be modified to~\( \abs{ \Im \rho } < T \) uniformly for~\( z \).

Now we compute the real part of the zero sum in~\eqref{ztSum}.
Suppose~\( y = \gamma_{0} \) for some zero~\( \rho_{0} = 1 / 2 + i \gamma_{0} \),
then there is a pair \( w_{x}\mleft( - \varepsilon \mright) \)~and~\(
\Re w_{x}\mleft( - \varepsilon - i \cdot 2 \gamma_{0} \mright) \).
The former detects this zero while the latter equals
\begin{equation*}
  \frac{1}{4 \gamma_{0}^{2} + \varepsilon^{2}}
  \mleft( \varepsilon - \varepsilon x^{- \varepsilon}
  \cos\mleft( 2 \gamma_{0} \log x \mright) + 2 \gamma_{0}
  x^{- \varepsilon} \sin\mleft( 2 \gamma_{0} \log x \mright) \mright)
  = O_{y}\mleft( \varepsilon + x^{- \varepsilon} \mright).
\end{equation*}
Now regardless of~\( y = \gamma_{0} \) or not, the rest of the sum is
\begin{equation}
  \label{zSum}
  \sum_{\substack{ 0 < \gamma < T \\ \gamma \ne \pm y }}
  \mleft( \frac{x^{iz + i \gamma} - 1}{i z + i \gamma}
    + \frac{x^{i z - i \gamma} - 1}{i z - i \gamma} \mright)
  = \sum_{\substack{ 0 < \gamma < T \\ \gamma \ne \pm y }}
  \frac{2 i}{z^{2} - \gamma^{2}}
  \mleft( z - x^{i z} \mleft( z \cos\mleft( \gamma \log x \mright)
  - i \gamma \sin\mleft( \gamma \log x \mright) \mright) \mright).
\end{equation}
The real part of the first term above is
\begin{equation}
  \label{leadBd}
  \Re \sum_{\substack{ 0 < \gamma < T \\ \gamma \ne \pm y }}
  \frac{2 i z}{z^{2} - \gamma^{2}}
  = 2 \varepsilon \sum_{\substack{ 0 < \gamma < T \\ \gamma \ne \pm y }}
  \frac{y^{2} + \gamma^{2} + \varepsilon^{2}}
  {\mleft( y^{2} - \gamma^{2} - \varepsilon^{2} \mright)^{2} + 4 \varepsilon^{2} y^{2}}
  \le 2 \varepsilon \sum_{0 < \gamma \ne \pm y}
  \frac{y^{2} + \gamma^{2} + 1}{\mleft( y^{2} - \gamma^{2} \mright)^{2}},
\end{equation}
while the cosine part in~\eqref{zSum} is
\begin{equation}
  \label{cosBd}
  \sum_{\substack{ 0 < \gamma < T \\ \gamma \ne \pm y }}
  \abs{ \frac{2 i z \cdot x^{iz} \cos\mleft( \gamma \log x \mright)}{z^{2} - \gamma^{2}} }
  \le 2 x^{- \varepsilon} \sum_{0 < \gamma \ne \pm y}
  \frac{\sqrt{ y^{2} + 1 }}{\abs{ y^{2} - \gamma^{2} }}.
\end{equation}
The only divergent part of~\eqref{zSum} is the sine,
which could be written as
\begin{equation}
  \label{sinBd}
  2 x^{iz} \sum_{\substack{ 0 < \gamma < T \\ \gamma \ne \pm y }}
  \mleft( \frac{1}{\gamma} - \frac{z^{2} / \gamma}{z^{2} - \gamma^{2}} \mright)
  \sin\mleft( \gamma \log x \mright)
  \ll x^{- \varepsilon} \mleft( \log^{2} T
    + \sum_{0 < \gamma \ne \pm y}
    \frac{y^{2} + 1}{\gamma \abs{ y^{2} - \gamma^{2} }} \mright),
\end{equation}
using the trivial bound
\begin{equation}
  \label{trivBd}
  \sum_{0 < \gamma < T} \frac{\abs{ \sin \mleft( \gamma \log x \mright) }}{\gamma}
  = \int_{0}^{T} \frac{1}{t} \df N\mleft( t \mright)
  \ll \log^{2} T.
\end{equation}
Together the real part of~\eqref{zSum} is
\begin{multline*}
  \Re \sum_{\substack{ 0 < \gamma < T \\ \gamma \ne \pm y }}
  \mleft( \frac{x^{iz + i \gamma} - 1}{i z + i \gamma}
  + \frac{x^{i z - i \gamma} - 1}{i z - i \gamma} \mright)
  \ll \varepsilon \sum_{0 < \gamma \ne \pm y}
  \frac{y^{2} + \gamma^{2} + 1}{\mleft( y^{2} - \gamma^{2} \mright)^{2}} \\
  + x^{- \varepsilon} \mleft( \sum_{0 < \gamma \ne \pm y}
    \frac{\sqrt{ y^{2} + 1 }}{\abs{ y^{2} - \gamma^{2} }}
    + \log^{2} T
    + \sum_{0 < \gamma \ne \pm y}
    \frac{y^{2} + 1}{\gamma \abs{ y^{2} - \gamma^{2} }} \mright).
\end{multline*}

Now let~\( T = x \log^{2} x \).
In the integral of~\eqref{ztSum},
the unnormalized sine cardinal picks out
half the rest of the integrand at~\( t = 0 \) as~\( x \to + \infty \)
(this step will not be uniform if~\( \varepsilon \to 1 / 2^{-} \)
due to the pole of digamma at the origin, while our~\( \varepsilon \to 0^{+} \) later).
Then the two digammas could be combined by replacing~\( i z \) by~\( - i y \)
at a cost of~\( O\mleft( \varepsilon \mright) \).
Together the real part of~\eqref{ztSum} reduces to
\begin{multline}
  \label{ztPf}
  \Re \mleft( \sum_{n \le x} \frac{\Lambda\mleft( n \mright)}{n^{1 / 2 + \varepsilon + i y}}
    - \frac{x^{1 / 2 - \varepsilon - i y} - 1}{1 / 2 - \varepsilon - i y}
    + \frac{x^{- 1 / 2 - \varepsilon - i y} - 1}{1 / 2 + \varepsilon + i y} \mright)
  - \frac{1}{2} \Re \frac{\Gamma'}{\Gamma}\mleft( \frac{1 / 2 - i y}{2} \mright)
  + \frac{1}{2} \log \pi \\
  =
  \begin{dcases}
    - w_{x}\mleft( - \varepsilon \mright) + O_{y}\mleft( \varepsilon \mright)
    + o_{y, \varepsilon}\mleft( 1 \mright)
    & \text{if}~y = \gamma, \\
    O_{y}\mleft( \varepsilon \mright)
    + o_{y, \varepsilon}\mleft( 1 \mright)
    & \text{otherwise}.
  \end{dcases}
\end{multline}
More explicitly the error is
\begin{equation}
  \label{err}
  \varepsilon \sum_{0 < \gamma \ne \pm y}
  \frac{y^{2} + \gamma^{2} + 1}{\mleft( y^{2} - \gamma^{2} \mright)^{2}}
  + x^{- \varepsilon}
  \mleft( \log^{2} x + \sum_{0 < \gamma \ne \pm y}
  \frac{y^{2} + 1}{\abs{ y^{2} - \gamma^{2} }} \mright)
  + \frac{1}{\log x}.
\end{equation}
By~\autoref{lmZ}, the first zero sum above is
\begin{equation}
  \label{zBd}
  \sum_{0 < \gamma \ne \pm y}
  \frac{y^{2} + \gamma^{2} + 1}{\mleft( y^{2} - \gamma^{2} \mright)^{2}}
  \le \mleft( 1 + \frac{1}{y^{2} + 1} \mright)
  \biggl( \, \sum_{\substack{ 0 < \gamma \ne \pm y \\ \abs{ \gamma - y } < 1 }}
    + \sum_{\substack{ 0 < \gamma \ne \pm y \\ \abs{ \gamma - y } \ge 1 }} \, \biggr)
  \frac{1}{\mleft( y - \gamma \mright)^{2}}
  \ll \frac{\log \mleft( 2 + \abs{ y } \mright)}{\min_{\gamma \ne y} \abs{ y - \gamma }^{2}},
\end{equation}
and the other sum could be similarly treated.
Now we choose~\( x^{\varepsilon} = \log^{3} x \) and
reduce the number of logarithms by the scale change~\( x \mapsto e^{x} \).
This completes the proof of~\autoref{tmZt}.

\subsection{Proof~of~\autoref{tmDr}}
A few changes are needed for Dirichlet~\( L \)-functions.
The polynomial~\eqref{dpDef} becomes
\begin{equation*}
  D_{x}\mleft( z, \chi \mright)
  \coloneqq \sum_{n \le x} \frac{\chi\mleft( n \mright) \Lambda\mleft( n \mright)}
  {n^{1 / 2 - i z}}
  = \frac{1}{2 \pi i} \int_{\mleft( c \mright)}
  - \frac{L'}{L}\mleft( 1 / 2 - i z + s, \chi \mright)
  \cdot w_{x}\mleft( s \mright) \df s.
\end{equation*}
Now~\( s = 1 \) is no longer a pole
while~\( s = 0 \) is an extra zero for even characters,
whence the residue formula~\eqref{rsd} becomes
\begin{equation}
  \label{dirRsd}
  I_{1} + I_{2} + I_{3} + I_{4}
  = - \left( 1 - \mathfrak{a}_{\chi} \right)
  w_{x}\mleft( - 1 / 2 + i z \mright)
  - \sum_{\abs{ \Im \rho_{\chi} - y } < T}
  w_{x}\mleft( \rho - 1 / 2 + i z \mright).
\end{equation}
The functional equation~\eqref{fneq} in this case is
\begin{equation*}
  - \frac{L'}{L}\mleft( 1 / 2 + s, \chi \mright)
  - \frac{L'}{L}\mleft( 1 / 2 - s, \overline{\chi} \mright)
  = \frac{1}{2} \frac{\Gamma'}{\Gamma}
  \mleft( \frac{1 / 2 + \mathfrak{a}_{\chi} + s}{2} \mright)
  + \frac{1}{2} \frac{\Gamma'}{\Gamma}
  \mleft( \frac{1 / 2 + \mathfrak{a}_{\chi} - s}{2} \mright)
  - \log \frac{\pi}{q}.
\end{equation*}
When moving the digamma integral to~\( \Re s = 0 \),
for even characters an extra pole exists at~\( s = - 1 / 2 + i z \),
canceling the extra zero in~\eqref{dirRsd}.
Then~\eqref{sinc} becomes
\begin{multline*}
  I_{3} = \frac{1}{2 \pi i} \int_{c - i T}^{c + i T}
  - \frac{L'}{L}\mleft( 1 / 2 + i z + s, \overline{\chi} \mright)
  \cdot w_{x}\mleft( - s \mright) \df s
  - \left( 1 - \mathfrak{a}_{\chi} \right) w_{x}\mleft( - 1 / 2 + i z \mright)
  + O\mleft( \frac{\log T}{T} \mright) \\
  - \frac{1}{2 \pi} \int_{- T}^{T}
  \mleft( \frac{1}{2} \frac{\Gamma'}{\Gamma}
    \mleft( \frac{1 / 2 + \mathfrak{a}_{\chi} + i z + it}{2} \mright)
  + \frac{1}{2} \frac{\Gamma'}{\Gamma}
  \mleft( \frac{1 / 2 + \mathfrak{a}_{\chi} - i z - it}{2} \mright)
  - \log \frac{\pi}{q} \mright) \cdot \frac{\sin\mleft( t \log x \mright)}{t} \df t,
\end{multline*}
while~\eqref{ztSum} becomes
\begin{multline}
  \label{chiSum}
  \sum_{n \le x} \frac{\chi\mleft( n \mright) \Lambda\mleft( n \mright)}{n^{1 / 2 - i z}}
  = - \sum_{\abs{ \Im \rho_{\chi} - y } < T}
  w_{x}\mleft( \rho_{\chi} - 1 / 2 + i z \mright)
  + O_{q}\mleft( \frac{x \log^{2} \mleft( x T \mright)}{T \log x} \mright) \\
  + \frac{1}{2 \pi} \int_{- T}^{T}
  \mleft( \frac{1}{2} \frac{\Gamma'}{\Gamma}
    \mleft( \frac{1 / 2 + \mathfrak{a}_{\chi} + i z + it}{2} \mright)
  + \frac{1}{2} \frac{\Gamma'}{\Gamma}
  \mleft( \frac{1 / 2 + \mathfrak{a}_{\chi} - i z - it}{2} \mright)
  - \log \frac{\pi}{q} \mright) \cdot \frac{\sin\mleft( t \log x \mright)}{t} \df t.
\end{multline}
Next we add above
to its conjugate version for~\( \overline{\chi} \)
at~\( \tilde{z} = y + i \varepsilon \), then apply
the zero symmetry~\( \rho_{\chi} \leftrightarrow 1 - \rho_{\overline{\chi}} \) to obtain
\begin{multline*}
  \Re \sum_{n \le x} \frac{\chi\mleft( n \mright) \Lambda\mleft( n \mright)}{n^{1 / 2 - i z}}
  = - \frac{1}{2} \sum_{\abs{ \Im \rho_{\chi} } < T}
  w_{x}\mleft( \rho_{\chi} - 1 / 2 + i z \mright)
  + w_{x}\mleft( 1 / 2 - \rho_{\chi} + i z \mright)
  + O_{q}\mleft( \frac{x \log^{2} \mleft( x T \mright)}{T \log x} \mright) \\
  + \frac{1}{2 \pi} \int_{- T}^{T}
  \mleft( \frac{1}{2} \frac{\Gamma'}{\Gamma}
    \mleft( \frac{1 / 2 + \mathfrak{a}_{\chi} + i z + it}{2} \mright)
  + \frac{1}{2} \frac{\Gamma'}{\Gamma}
  \mleft( \frac{1 / 2 + \mathfrak{a}_{\chi} - i z - it}{2} \mright)
  - \log \frac{\pi}{q} \mright) \cdot \frac{\sin\mleft( t \log x \mright)}{t} \df t.
\end{multline*}

If~\( L\mleft( s, \chi \mright) \)
vanishes at~\( s = 1 / 2 \), when~\( y = 0 \)
there is only a single term~\( w_{x}\mleft( - \varepsilon \mright) \).
Moreover in~\eqref{zSum} there is an extra term
\begin{equation*}
  \Re w_{x}\mleft( - \varepsilon - i y \mright)
  = \frac{1}{y^{2} + \varepsilon^{2}}
  \mleft( \varepsilon - \varepsilon x^{- \varepsilon}
  \cos\mleft( y \log x \mright)
  + y x^{- \varepsilon} \sin\mleft( y \log x \mright) \mright)
  = O\mleft( \varepsilon y^{-2}
  + x^{- \varepsilon} \log x \mright),
\end{equation*}
which could be inserted into the error~\eqref{err}
by extending the first zero sum to include the case~\( \Im \rho_{\chi} = 0 \).
Finally, the formula corresponding to~\eqref{ztPf} is
\begin{multline}
  \label{drPf}
  \Re \sum_{n \le x} \frac{\chi\mleft( n \mright) \Lambda\mleft( n \mright)}
  {n^{1 / 2 + \varepsilon + i y}}
  - \frac{1}{2} \Re \frac{\Gamma'}{\Gamma}\mleft( \frac{1 / 2 + \mathfrak{a}_{\chi} - i y}{2} \mright)
  + \frac{1}{2} \log \frac{\pi}{q} \\
  =
  \begin{dcases}
    - w_{x}\mleft( - \varepsilon \mright)
    + O_{q, y}\mleft( \varepsilon \mright)
    + o_{q, y, \varepsilon}\mleft( 1 \mright)
    & \text{if}~y = \gamma_{\chi}, \\
    O_{q, y}\mleft( \varepsilon \mright)
    + o_{q, y, \varepsilon}\mleft( 1 \mright)
    & \text{otherwise.}
  \end{dcases}
\end{multline}
By~\autoref{lmC}, the zero sum bound~\eqref{zBd} becomes
\begin{equation*}
  \sum_{0 \le \gamma_{\chi} \ne \pm y}
  \frac{y^{2} + \gamma_{\chi}^{2} + 1}{\mleft( y^{2} - \gamma_{\chi}^{2} \mright)^{2}}
  \ll_{q} \frac{\log \mleft( 2 + \abs{ y } \mright)}
  {\min_{\gamma_{\chi} \ne y} \abs{ y - \gamma_{\chi} }^{2}},
\end{equation*}
and this completes the proof of~\autoref{tmDr}.

\section{Decomposition}

The notation and labeling in this section
have different meaning for different families of~\( L \)-functions under discussion.
\autoref{crInd} follows from considering \autoref{tmZt}~and~\autoref{tmDr} together.
Under~\( 0 \)-indexing,
let~\( \mleft\{ \chi_{j} \mright\} \) run through the~\( \varphi\mleft( q \mright) \)
Dirichlet characters with~\( \chi_{0} \) principal,
and~\( \mleft\{ a_{k} \mright\} \) run through~\( \mleft( \mathbb{Z} / q \mathbb{Z} \mright)^{\times} \)
with~\( a_{0} = 1 \). For~\( 0 \le j, k \le \varphi\mleft( q \mright) - 1 \) let
\begin{equation}
  \label{cMat}
  \mleft( M \mright)_{j, k}
  = \chi_{j}\mleft( a_{k} \mright)
  ~\text{and}~
  \mleft( \bm{S} \mright)_{k}
  = \sum_{\substack{ 0 < n \log p \le x \\ p^{n} \equiv a_{k} \bmod q }}
  \frac{\log p}{p^{n \mleft( 1 / 2 + \varepsilon_{x} + i y \mright)}}.
\end{equation}
Then the indicators \eqref{indZt}~and~\eqref{indDr}
together give the linear system
\begin{equation}
  \label{ls}
  \Re \mleft( M \cdot \bm{S}^{\intercal} \mright)
  =
  - \begin{pmatrix}
    \mathbb{1}_{0} \\ \mathbb{1}_{1} \\
    \vdots \\ \mathbb{1}_{\varphi\mleft( q \mright) - 1}
  \end{pmatrix}
  + \Re N
  \begin{pmatrix}
    1 \\ 0 \\
    \vdots \\ 0
  \end{pmatrix}
  + G_{\textnormal{e}}
  \begin{pmatrix}
    1 \\ \vdots
  \end{pmatrix}
  + G_{\textnormal{o}}
  \begin{pmatrix}
    0 \\ \vdots
  \end{pmatrix}
  - \frac{1}{2} \log \frac{\pi}{q} \cdot \bm{1}^{\intercal},
\end{equation}
where~\( N \) is the exponential noise present only in~\eqref{indZt}
and \( G_{\textnormal{e}}, G_{\textnormal{o}} \) the digammas,
\begin{align*}
  N
  & = \frac{x^{1 / 2 - \varepsilon_{x} - i y} - 1}
  {1 / 2 - \varepsilon_{x} - i y}
  - \frac{x^{- 1 / 2 - \varepsilon_{x} - i y} - 1}
  {1 / 2 + \varepsilon_{x} + i y}, \\
  G_{\textnormal{e}}
  & = \frac{1}{2} \Re \frac{\Gamma'}{\Gamma}
  \mleft( \frac{1 / 2 - i y}{2} \mright)
  ~\text{and}~G_{\textnormal{o}}
  = \frac{1}{2} \Re \frac{\Gamma'}{\Gamma}
  \mleft( \frac{3 / 2 - i y}{2} \mright).
\end{align*}
Note different orders of \( \mleft\{ \chi_{j} \mright\} \)~and~\( \mleft\{ a_{k} \mright\} \)
simply permute the characters and residue classes.

For a primitive character~\( \chi \) on the~\( j \)-th row,
the equation is~\eqref{indDr} where the spike function
on the right is denoted by \( \mathbb{1}_{\chi} \)~or~\( \mathbb{1}_{j} \),
while ideally~\( \mathbb{1}_{\chi}\mleft( y \mright)
= \sum_{\gamma_{\chi}} \delta\mleft( y - \gamma_{\chi} \mright) \).

For a non-primitive character~\( \chi' \) on the~\( j \)-th row induced by some~\( \chi \),
its indicator~\( \mathbb{1}_{\chi'} \) is no longer pure spikes,
but rather~\( \mathbb{1}_{\chi} \) plus some elementary factors.
More explicitly when~\( q = p^{k} \), the series difference
between the principal character
and the Riemann zeta version~\eqref{indZt} are
powers of~\( p \) and the constant~\( \tfrac{1}{2} \log p^{k} \),
giving the compensation
\begin{equation}
  \label{pCor}
  \begin{aligned}
    \mathbb{1}_{0}\mleft( y \mright)
    & = \mathbb{1}_{\zeta}\mleft( y \mright)
    + \frac{k}{2} \log p
    + \Re \sum_{1 \le n \le \mleft\lfloor \frac{\log x}{\log p} \mright\rfloor}
    \frac{\log p}{\mleft( p^{n} \mright)^{1 / 2 + \varepsilon_{x} + i y}} \\
    & = \mathbb{1}_{\zeta}\mleft( y \mright)
    + \frac{\log p}{2}
    \Re \mleft( \frac{1 + p^{- 1 / 2 - \varepsilon_{x} - i y}}
    {1 - p^{- 1 / 2 - \varepsilon_{x} - i y}} \mright)
    + \frac{k - 1}{2} \log p
    + O\mleft( \frac{\log p}{x^{1 / 2 + \varepsilon_{x}}} \mright).
  \end{aligned}
\end{equation}
For a non-principal non-primitive character~\( \chi' \bmod p^{k} \)
induced by~\( \chi \bmod p^{\ell} \),
the series remains the same since~\( \mleft( n, p^{k} \mright) = 1 \)
is equivalent to~\( \mleft( n, p^{\ell} \mright) = 1 \),
so only the constant needs compensation, giving
\begin{equation}
  \label{npCor}
  \mathbb{1}_{\chi'}\mleft( y \mright)
  = \mathbb{1}_{\chi}\mleft( y \mright) + \frac{k - \ell}{2} \log p.
\end{equation}
For a general composite modulus,
the compensation needed is a combination of \eqref{pCor}~and~\eqref{npCor}.
Let~\( \chi' \bmod q' \) be induced by~\( \chi \bmod q \), then
\begin{multline}
  \label{chiCor}
  \mathbb{1}_{\chi'}\mleft( y \mright)
  = \mathbb{1}_{\chi}\mleft( y \mright)
  + \sum_{\substack{ p \mid q' \\ p \nmid q }} \frac{\log p}{2}
  \Re \mleft( \frac{1 + p^{- 1 / 2 - \varepsilon_{x} - i y}}
  {1 - p^{- 1 / 2 - \varepsilon_{x} - i y}} \mright)
  + \sum_{\substack{ p^{k} \parallel q' \\ p \nmid q }} \frac{k - 1}{2} \log p \\
  + \sum_{\substack{ p^{k} \parallel q' \\ p^{\ell} \parallel q }}
  \frac{k - \ell}{2} \log p
  + O\mleft( \frac{\log q}{x^{1 / 2 + \varepsilon_{x}}} \mright),
\end{multline}
where~\( p^{k} \parallel q \) means~\( p^{k} \mid q \) but~\( p^{k + 1} \nmid q \).
In particular for a principal~\( \chi \),
\begin{equation}
  \label{ztCor}
  \mathbb{1}_{0}\mleft( y \mright)
  = \mathbb{1}_{\zeta}\mleft( y \mright)
  + \sum_{p \mid q} \frac{\log p}{2}
  \Re \mleft( \frac{1 + p^{- 1 / 2 - \varepsilon_{x} - i y}}
  {1 - p^{- 1 / 2 - \varepsilon_{x} - i y}} \mright)
  + \sum_{p^{k} \parallel q} \frac{k - 1}{2} \log p
  + O\mleft( \frac{\log q}{x^{1 / 2 + \varepsilon_{x}}} \mright).
\end{equation}

The matrix \( M \) is invertible as~\( M \cdot M^{*}
= \varphi\mleft( q \mright) \cdot I \) by character orthogonality,
with~\( M^{*} \) its conjugate transpose.
\autoref{crInd} follows from~\( \mleft( \bm{S} \mright)_{0} \)
by inverting the matrix equation~\eqref{ls},
while the extra error~\( O_{q}\mleft( 1 \mright) \)
contains finitely many elementary factors in \eqref{chiCor}~and~\eqref{ztCor}
from non-primitive characters.
Note this step is only complete for rows corresponding to
real characters, because for complex ones the imaginary part of~\eqref{ls}
is needed, see~\autoref{secIm}.
This inversion also demonstrates how components of
explicit formulae are distributed among residue classes.
For example, the pole~\( s = 1 \) of~\( \zeta\mleft( s \mright) \)
gives Dirichlet's theorem, while a Siegel zero distorts
this by almost doubling and halving them alternatively.

The geometric series in~\eqref{pCor} and its variants
have appeared in the study of zero periodicity
by~\citeauthor{fordDistributionImaginaryParts2009}
for the Riemann zeta function~\cite{fordDistributionImaginaryParts2005}
and Selberg class~\cite{fordDistributionImaginaryParts2009}.
Here with the missing~\( 0 \)-component it is completed to a full Poisson kernel,
which could be regarded as the contribution from a single prime to the zeros.

More abstractly speaking, the columns of a DFT matrix are eigenvectors of
corresponding permutation or circulant matrices,
while eigenvalues are roots of unity or their linear combinations.
This means that for any~\( q \in \mathbb{N} \),
there are operators acting on Dirichlet series
by permutation or circular shift within each prime component
of the~\( \varphi\mleft( q \mright) \) residue classes,
and the sums in~\eqref{cMat} or their linear combinations are the eigenfunctions.

In addition to the full decompositon described above
there are also sub-decompositions, corresponding to intermediate fields
between~\( \mathbb{Q} \) and some cyclotomic field.
For example, we pick a primitive and quadratic~\( \chi \bmod q \) and let
\begin{equation*}
  S_{\pm} =
  \sum_{\substack{ 0 < n \log p \le x \\ \chi\mleft( p^{n} \mright) = \pm 1 }}
  \frac{\log p}{p^{n \mleft( 1 / 2 + \varepsilon_{x} + i y \mright)}}.
\end{equation*}
The linear system~\eqref{ls} in this case is
\begin{equation}
  \label{qd}
  \Re
  \begin{pmatrix}
    1 & 1 \\
    1 & -1
  \end{pmatrix}
  \begin{pmatrix}
    S_{+} \\ S_{-}
  \end{pmatrix}
  =
  -
  \begin{pmatrix}
    \mathbb{1}_{0} \\ \mathbb{1}_{\chi}
  \end{pmatrix}
  + \Re N
  \begin{pmatrix}
    1 \\ 0
  \end{pmatrix}
  + G_{\textnormal{e}}
  \begin{pmatrix}
    1 \\ \vdots
  \end{pmatrix}
  + G_{\textnormal{o}}
  \begin{pmatrix}
    0 \\ \vdots
  \end{pmatrix}
  - \frac{1}{2} \log \frac{\pi}{q}
  \begin{pmatrix}
    1 \\ 1
  \end{pmatrix},
\end{equation}
where~\( \mathbb{1}_{0} \) is given by~\eqref{ztCor}
and \( \chi\mleft( -1 \mright) \) determines
taking \( G_{\textnormal{e}} \)~or~\( G_{\textnormal{o}} \).

As another example, let~\( \mleft( q_{1}, q_{2} \mright) = 1 \),
we fix a primitive~\( \chi \bmod q_{1} \)
and let~\( \mleft\{ \chi_{j} \mright\} \) run through the characters modulo~\( q_{2} \)
with~\( \chi_{0} \) principal, and~\( \mleft\{ a_{k} \mright\} \) run
through~\( \mleft( \mathbb{Z} / q_{2} \mathbb{Z} \mright)^{\times} \) with~\( a_{0} = 1 \).
Let~\( M \) be given by~\eqref{cMat} for~\( q = q_{2} \) and let
\begin{equation*}
  \mleft( \bm{S} \mright)_{k} = \sum_{\substack{ 0 < n \log p \le x \\ p^{n} \equiv a_{k} \bmod q_{2} }}
  \frac{\chi_{1}\mleft( p^{n} \mright) \log p}
  {p^{n \mleft( 1 / 2 + \varepsilon_{x} + i y \mright)}}.
\end{equation*}
Then in the full linear system~\eqref{ls},
rows corresponding to characters containing~\( \chi \) as a factor
could be collected to form the sub-system
\begin{equation*}
  \Re \mleft( M \cdot \bm{S}^{\intercal} \mright)
  =
  - \begin{pmatrix}
    \mathbb{1}_{0} \\
    \vdots \\ \mathbb{1}_{\varphi\mleft( q_{2} \mright) - 1}
  \end{pmatrix}
  + G_{\textnormal{e}}
  \begin{pmatrix}
    1 \\ \vdots
  \end{pmatrix}
  + G_{\textnormal{o}}
  \begin{pmatrix}
    0 \\ \vdots
  \end{pmatrix}
  - \frac{1}{2} \log \frac{\pi}{q_{1} q_{2}} \cdot \bm{1}^{\intercal},
\end{equation*}
where~\( \mathbb{1}_{k} \) for~\( 0 \le k \le \varphi\mleft( q_{2} \mright) - 1 \)
are compensated similarly as in~\eqref{chiCor}
by comparing series. The noise~\( N \) is absent
because the DC component is a Dirichlet~\( L \)-function now.

\section{Imaginary Part}
\label{secIm}

Throughout the proof we have kept complex notation,
so it is natural to attempt at the imaginary part
of the zero sum~\eqref{zSum} as well.
Note the real part of digamma is inherent by~\eqref{ztSum},
and the absolute value bounds \eqref{cosBd}~and~\eqref{sinBd} for oscillatory terms
could be reused. Together with~\( \Im w_{x}\mleft( - \varepsilon \mright) = 0 \) and
\begin{align*}
  \Im w_{x}\mleft( - \varepsilon - 2 i \gamma_{0} \mright)
  & = \frac{- 2 \gamma_{0} + x^{- \varepsilon}
  \mleft( 2 \gamma_{0} \cos\mleft( 2 \gamma_{0} \log x \mright)
  + \varepsilon \sin\mleft( 2 \gamma_{0} \log x \mright) \mright)}
  {4 \gamma_{0}^{2} + \varepsilon^{2}}
  = O_{y}\mleft( 1 \mright), \\
  \Im \sum_{\substack{ 0 < \gamma < T \\ \gamma \ne \pm y }}
  \frac{2 i z}{z^{2} - \gamma^{2}}
  & = - 2 y \sum_{\substack{ 0 < \gamma < T \\ \gamma \ne \pm y }}
  \frac{y^{2} - \gamma^{2} + \varepsilon^{2}}
  {\mleft( y^{2} - \gamma^{2} - \varepsilon^{2} \mright)^{2} + 4 \varepsilon^{2} y^{2}}
  = O_{y}\mleft( 1 \mright),
\end{align*}
we find the imaginary part is~\( O_{y}\mleft( 1 \mright) \)
no matter at zeros or bounded away, not too illuminating.
To see its actual behavior, the following limit case is more appropriate,
\begin{equation}
  \label{lmt}
  \lim_{\varepsilon \to 0^{+}}\frac{x^{- \varepsilon + i \mleft( y - \gamma \mright)} - 1}
  {- \varepsilon + i \mleft( y - \gamma \mright)}
  =
  \begin{dcases}
    \log x + i \cdot 0 & \text{if}~y = \gamma, \\
    \frac{\sin\mleft( \mleft( y - \gamma \mright) \log x \mright)}{y - \gamma}
    + i \cdot \frac{1 - \cos\mleft( \mleft( y - \gamma \mright) \log x \mright)}{y - \gamma}
    & \text{otherwise}.
  \end{dcases}
\end{equation}
At a zero~\( \gamma \) the real part is~\( \log x \),
whereas the imaginary part vanishes but attains approximately~\( \pm 0.72 \log x \)
at~\( \gamma \pm 2.33 / \log x \), overshooting and undershooting
in the same manner as the Wilbraham--Gibbs phenomenon.
A similar analysis shows that Landau's result~\eqref{ld},
while technically correct, actually needs a real part on the left.
We do not proceed any further since
the current proof fails to give the correct size even for the real part.
Instead see~\autoref{fgReIm} for a comparison of their behaviors.

Under RH \citeauthor{kharePhaseRiemannZeta1997}~\cite[\S 5]{kharePhaseRiemannZeta1997}
showed that~\( \Re \zeta'/\zeta\mleft( 1 / 2 + i t \mright) \)
approximates~\( \sum_{\gamma} \delta\mleft( t - \gamma \mright) \)
plus some elementary factors.
Together with~\eqref{lmt}, they suggest what happens on the critical line
is simply \citeauthor{whittakerXVIIIonFunctionsWhich1915}--\citeauthor{shannonCommunicationPresenceNoise1949}
interpolation at zeros.
In fact, the real and imaginary parts are related via Hilbert transform
by the Sokhotski--Plemelj theorem on~\( \mathbb{R} \),
also known as the Kramers--Kronig relation,
\begin{equation}
  \label{kk}
  \frac{1}{\pi i} \lim_{\varepsilon \to 0^{+}}
  \frac{1}{\pm \varepsilon + i \mleft( y - \gamma \mright)}
  = \pm \delta\mleft( y - \gamma \mright)
  - i \cdot \frac{1}{\pi}\pv\mleft( \frac{1}{y - \gamma} \mright),
\end{equation}
where \( \pv \) denotes the Cauchy principal value.
When the linear system~\eqref{ls} is interpreted ideally in the complex sense,
the~\( \mathbb{1}_{k} \)'s should be replaced by the corresponding sums of~\eqref{kk}.

\section{Numerical Experiment}

In this section we drop the offset~\( \varepsilon_{x} \),
as it requires all prime powers up to~\( 8.8 \times 10^{19} \)
for~\( \varepsilon_{x} \le 1 / 4 \). Instead we choose the more reasonable
cutoff~\( e^{x} = 10^{9} \) whence the predicted spike size is~\( x \approx 20.72 \).
We provide uniform versions where the cutoff is fixed while the ordinate~\( y \) varies.
In some cases we also give pointwise versions where~\( y \) is fixed
while the cutoff varies from~\( 2 \) to~\( 10^{9} \) in logarithmic scale.
Pictures are made with~\texttt{Julia}~\cite{bezansonJuliaFreshApproach2017}
using \emph{Fast~Fourier~Transform}~(FFT) algorithms
in \citeauthor{cooleyAlgorithmMachineCalculation1965}~\cite{cooleyAlgorithmMachineCalculation1965}
as well as \citeauthor{goodInteractionAlgorithmPractical1958}~\cite{goodInteractionAlgorithmPractical1958}.

\begin{figure}[htbp]
  \centering
  \includegraphics[width=\linewidth]{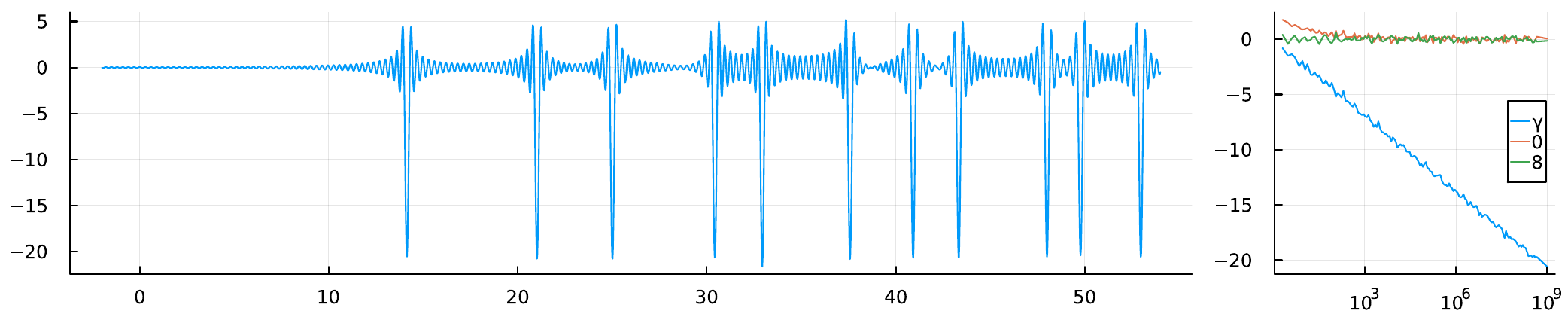}
  \caption{Zero indicator of~\( \zeta\mleft( s \mright) \) around origin.}
  \label{fgZtL}
\end{figure}

In~\autoref{fgZtL}
we plot~\autoref{tmZt} for~\( y \in \mleft[ -2, 54 \mright] \),
and pick the first zero~\( \gamma = 14.13 \ldots \),
the origin and~\( y = 8 \).
The flat line around~\( y = 0 \) suggests the sum
\begin{equation*}
  \sum_{0 < n \log p \le x} \frac{\log p}{p^{n / 2}}
  - 4 \sinh\mleft( \frac{x}{2} \mright)
  - \frac{1}{2} \Re \frac{\Gamma'}{\Gamma}\mleft( \frac{1}{4} \mright)
  + \frac{\log \pi}{2}
  = - \sum_{\gamma} \frac{\sin\mleft( \gamma x \mright)}{\gamma}
\end{equation*}
is actually much smaller than the trivial bound~\( O\mleft( x^{2} \mright) \).
The formula above follows from~\eqref{ztSum} at~\( z = 0 \),
see also~\citeauthor{guinandSummationFormulaTheory1948}~\cite[(2.7)]{guinandSummationFormulaTheory1948}.

\begin{figure}[htbp]
  \centering
  \includegraphics[width=\linewidth]{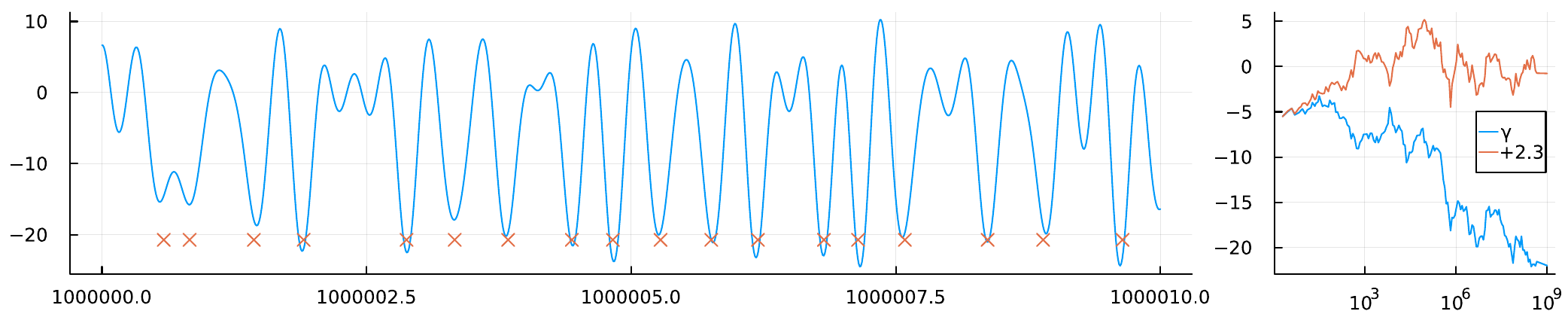}
  \caption{Zero indicator of~\( \zeta\mleft( s \mright) \) higher up.}
  \label{fgZtH}
\end{figure}

In~\autoref{fgZtH} we plot~\autoref{tmZt}
for~\( y \in \mleft[ 10^{6}, 10^{6} + 10 \mright] \),
and pick \( \gamma = 10^{6} + 1.90 \ldots \)~and~\( y = 10^{6} + 2.3 \).
Actual zeros are marked with crosses at predicted heights.

\begin{figure}[htbp]
  \centering
  \includegraphics[width=\linewidth]{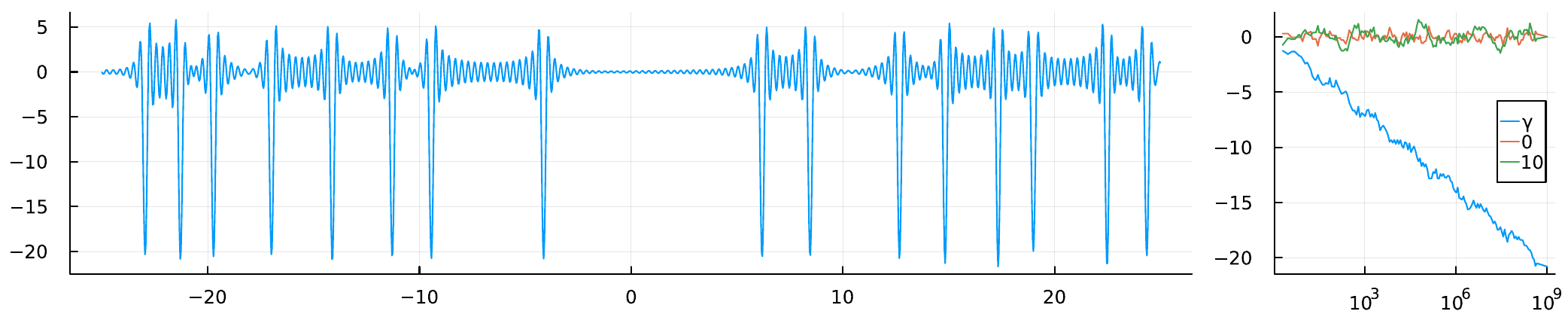}
  \caption{Zero indicator of~\( L\mleft( s, \chi \mright) \).}
  \label{fgDr}
\end{figure}

In~\autoref{fgDr} we choose the complex character~\( \chi \bmod 5 \)
determined by~\( \chi\mleft( 2 \mright) = i \)
and plot~\autoref{tmDr} for~\( y \in \mleft[ - 25, 25 \mright] \).
The ordinates picked are the first positive zero~\( \gamma_{\chi} = 6.18 \ldots \),
the origin and~\( y = 10 \).

The examples above are continuations of the work
by~\citeauthor{mumfordNumbersWorldEssays2023}~\cite[\S 11]{mumfordNumbersWorldEssays2023},
originally published as a blog
post\footnote{\url{https://www.dam.brown.edu/people/mumford/BookBlogPosts/Rhythms\%20of\%20the\%20Primes.pdf}}.

\begin{figure}[htbp]
  \centering
  \includegraphics[width=\linewidth]{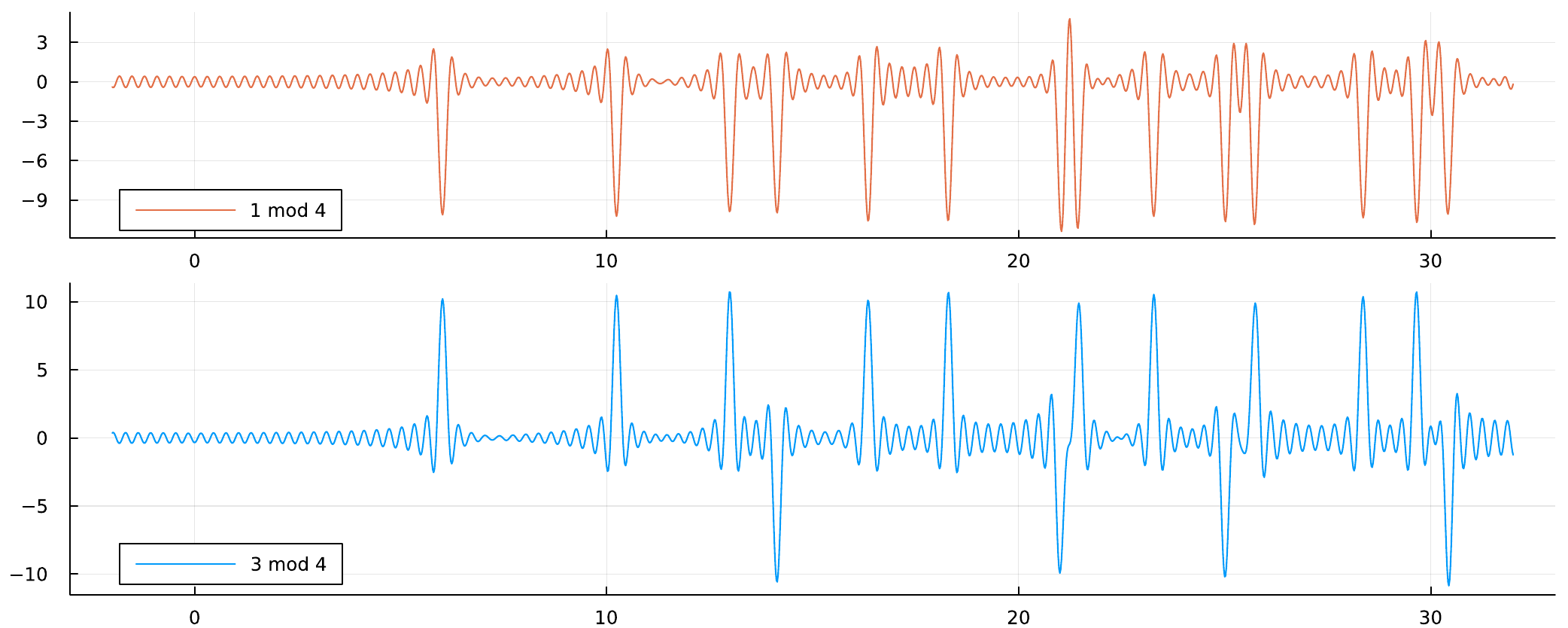}
  \caption{Quadratic decomposition.}
  \label{fgPm}
\end{figure}

Next in~\autoref{fgPm}
we choose~\( q = 4 \) and plot for~\( y \in \mleft[ -2, 32 \mright] \),
giving a simultaneous example for~\autoref{crInd}
and the quadratic decomposition~\eqref{qd}.
The two residue classes~\( \pm 1 \bmod 4 \)
give~\( - \mathbb{1}_{\zeta} \mp \mathbb{1}_{\chi} \) respectively,
with halved spike size~\( x / 2 \approx 10.36 \).
Exact compensations in \eqref{chiCor}~and~\eqref{ztCor}
are used to eliminate the extra error~\( O_{q}\mleft( 1 \mright) \).

\begin{figure}[htbp]
  \centering
  \includegraphics[width=\linewidth]{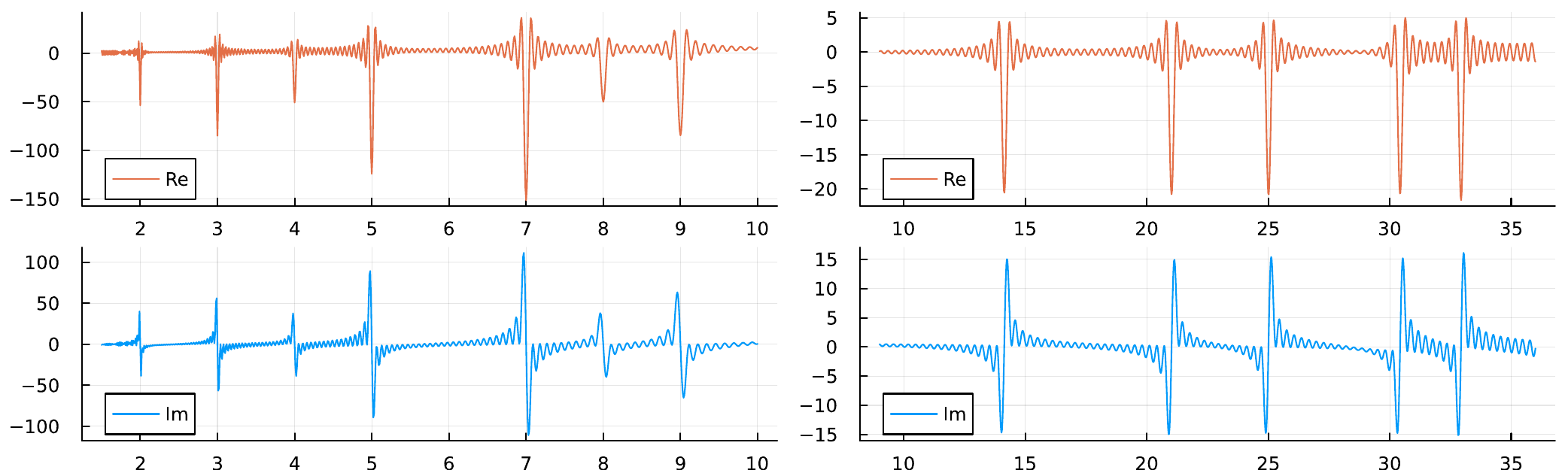}
  \caption{Real vs. Imaginary parts of \eqref{ld}~and~\eqref{indZt}.}
  \label{fgReIm}
\end{figure}

Lastly in~\autoref{fgReIm} we plot real and imaginary parts of
Landau's result~\eqref{ld} with~\( T = 500 \)
and~\eqref{indZt}. We choose a small~\( T \) to show the convergence
because in this direction it is linear rather than logarithmic.
See~\autoref{secIm} for discussions.

\printbibliography

\end{document}